\documentclass[b5,draft]{article}
\usepackage{amssymb,amsmath,amsfonts,curves,ascmac,fancybox,makeidx,color,amsthm,amsxtra,bm}
\allowdisplaybreaks

\def\qed{\hfill $\Box$}

\newtheorem{thm}{Theorem}       
\newtheorem{lem}{Lemma}         

\title{Behaviour of the least root modulo a prime of a polynomial}
\author{Yoshiyuki Kitaoka
\thanks{
\tt kitaoka@meijo-u.ac.jp
} 
} 
\date{}
\begin{document}
\maketitle 
Let 
\begin{equation}\label{eq0}
f(x)=x^n+a_{n-1}x^{n-1}+\dots+a_0
\end{equation} 
be  a monic  irreducible polynomial of degree $n\,(>1)$ with $a_i\in \mathbb Z$.
We put
\begin{equation*}
Spl_X(f):=\left\{p\le X\mid f(x) \text{ is fully splitting modulo
}p\right\}
\end{equation*}
for a positive number $X$ and $Spl(f):=Spl_\infty(f)$.
Here a letter $p$ denotes a prime number.
We require the following conditions on the local roots
$r_1,\dots,r_n\,( \in \mathbb{Z})$
 of $f(x)\equiv 0\bmod p$ for   a prime $p\in Spl(f)$\,:
\begin{align}
\label{eq1}
&f(x) \equiv \prod_{i=1}^n(x-r_i) \bmod p,
\\
\label{eq2}
&0\le r_1\le r_2\le\dots\le r_n<p.
\end{align}
The condition \eqref{eq2} determines  local roots $r_i$ uniquely.
In this note, we suppose that a polynomial $f$ has only a trivial linear relation among roots, that is
 for  complex roots $\alpha_1,\dots,\alpha_n$  of a polynomial $f$ in \eqref{eq0}
 a linear relation
 \begin{equation}\label{eq3}
\sum_{i=1}^n m_i\alpha_i=m\quad(m_i,m\in\mathbb{Q}).
\end{equation}
is satisfied only if $m_1=\dots=m_n$.
We know that if for an irreducible polynomial $f$  the degree $n$  is prime or the Galois group is isomorphic to 
the symmetric group $S_n$,
then there is only a trivial linear relation among roots \cite{K2}.

We conjectured  a kind of uniformity on the distribution of $(r_1/p,\dots,r_n/p)\linebreak[3]
\in[0,1)^n$ $(p\in Spl(f))$ in \cite{K1},  \cite{K2}.
Put
\begin{equation}\label{eq4}
\hat{\mathfrak{D}}_n:=
\{
(x_1\dots,x_n)\in [0,1)^n\mid
0\le x_1\le\dots\le x_n<1,\sum_{i=1}^n x_i\in\mathbb{Z}
\}
\end{equation}
and for a domain $D\subset[0,1)^n$ with $D=\overline{D^\circ}$
\begin{equation*}
Pr_D(f,X):=\displaystyle\frac{\#\{p\in Spl_X(f)\mid (r_1/p,\dots,r_n/p)\in D\}}{\#Spl_X(f)},
\end{equation*}
where local roots $r_i$ satisfy properties \eqref{eq1}, \eqref{eq2}.

\noindent
The conjecture is : Under an assumption that  a polynomial $f$  
has only a trivial linear relation \eqref{eq3} among roots,
\begin{equation}\label{eq5}
Pr_D(f):=\lim_{X\to\infty}Pr_D(f,X)
=\displaystyle\frac{ vol(  {D}\cap{\hat{\mathfrak{D}}_n})} {vol(\hat{\mathfrak{D}}_n)}.
\end{equation}
Here we give the right-hand side of \eqref{eq5} explicitly for
\begin{align*}
&D^a:=\{(x_1,\dots,x_n)\in [0,1)^n\mid a \le x_1 \} \quad(0\le a <1),
\\
&D_a:=\{(x_1,\dots,x_n)\in [0,1)^n\mid  x_1 \le a\} \quad(0\le a <1).
\end{align*}
We note that $Pr_{D_a}(f)$ is the density of primes $p$ satisfying $r_1/p<a$,
and 
$$
vol({D_a}\cap{\hat{\mathfrak{D}}_n})+  vol({D^a}\cap{\hat{\mathfrak{D}}_n})=
 vol(\hat{\mathfrak{D}}_n). 
 $$
\begin{thm}
The right-hand side $V^a$ of \eqref{eq5}  for $D=D^a$ is equal to
 $$
 \frac{1}{(n-1)!}         \sum_{2\le i\le n\atop \max{(ia,i-n+1)}\le k \le i-1}  (-1)^{n+i}{\binom{n}i}
(k-ia)^{n-1} .
 $$
 \end{thm}
This means that the conjectural value  of the density of primes $p$ satisfying that
every root $r\,(0\le r<p)$ of $f(x)\equiv0\bmod p$ is greater than $ap$ is equal to  $V^a$,
if a polynomial $f$ has only a trivial linear relation among roots.

Let us give more explicitly an expected density $V_a:= 1-V^a$ for a domain $D_a$ in the case of $n=2,3,4$.

\noindent
In case of $n=2$:
$$
\begin{cases}
2a&(0\le a \le1/2),
\\
1&(1/2\le a\le 1).
\end{cases}
$$
In case of $n=3$:
$$
\begin{cases}
3a-3a^2&(0\le a \le1/3),
\\
(1+3a^2)/2&(1/3\le a\le 1/2),
\\
-1+6a-9a^2/2&(1/2\le a\le2/3),
\\
1&(2/3\le a \le1).
\end{cases}
$$
In case of $n=4$:
$$
\begin{cases}
4a^3 - 6a^2 + 4a       & (0 \le a \le1/4),\\
-20a^3/3  + 2a^2 + 2a + 1/6  & (1/4 < a \le 1/3),\\
34a^3/3 - 16a^2 + 8a - 1/2    & (1/3 < a \le 1/2),\\
-22a^3/3 + 12a^2 - 6a + 11/6   & (1/2 < a \le 2/3),\\
32a^3/3 - 24a^2 + 18a - 7/2  & (2/3 < a \le3/4),\\
1  & (3/4 < a < 1).
\end{cases}
$$
Let us give  numerical data of differences between the  expected density $V^a$ and $Pr_{D^a}(f,X)$
$(a=1/(10n),\dots,(10n-1)/(10n))$ for $f = x^n+3x+1$,
whose Galois group is isometric to the symmetric group.
The following is a table\footnote{ Data were made by pari/gp. The PARI~Group, PARI/GP version
{\tt 2.8.0}, Bordeaux, 2014, http://pari.math.u-bordeaux.fr/.} of $\max_{1\le k \le10n-1} | Pr_{D^{k/(10n)}}(f,X_m) - V^{k/(10n)}|$,
where an integer $X_m$ denotes the least prime number exceeding $10^m$ in $Spl(f)$.
$$
\begin{array}{|c|c|c|c|c|c|c|c|c|c|}
\hline
&n=2&3&4&5&6&7&8&9&10\\
\hline
m=7&1(3)&1(3)&2(3)&9(3)&1(2)&5(2)&2(1)&2(1)&8(1)
\\
\hline
8&4(4)&6(4)&1(3)&6(3)&1(2)&2(2)&3(2)&1(1)&3(1)
\\
\hline
9&4(5)&3(4)&3(4)&1(3)&1(3)&8(3)&1(2)&3(2)&1(1)
\\
\hline
\end{array}
$$
In the table, $1(3)$ for $m=7,n=2$ means  $1.\cdots\times10^{-3}=0.001\cdots$ and so on.

Before a calculation of the volume of $vol(  {D^a}\cap{\hat{\mathfrak{D}}_n})$,
we refer to the following fundamental lemma (\cite{Fe}).
\begin{lem}
 For a natural number $n$, the volume of a subset of the unit cube $[0,1)^{n}$ 
 defined by $\{ (x_1,\dots,x_n)\in[0,1)^n \mid \,x_1+\dots+x_n\le x\}$ is given by
\begin{equation*}
 U_n(x):=\frac{1}{n!}\sum_{i=0}^n (-1)^i{n\choose i}\max(x-i,0)^n.
\end{equation*}
\end{lem}
Denoting by $\theta$ the angle of two hyperplanes of $\mathbb{R}^n$ defined by $x_1+\dots+x_n=0$ and $x_n=0$,
we have
\begin{align*}
&\hspace{5mm}vol({D^a}\cap{\hat{\mathfrak{D}}_n})
\\
&=vol(\{(x_1,\dots,x_n)\mid a\le x_1\le\dots <x_n<1,\sum x_i\in\mathbb{Z}\})
\\
&=\frac{1}{n!}vol(\{(x_1,\dots,x_n)\mid a\le x_1,\dots,x_n<1,\sum x_i\in\mathbb{Z}\})
\\
\intertext
{by $\sum_{i=1}^n x_i=\lceil \sum_{i=1}^{n-1} x_i\rceil$ and $x_n=\lceil\sum_{i=1}^{n-1} x_i\rceil-\sum_{i=1}^{n-1} x_i$, where $k:=\lceil x\rceil$ is an  integer satisfying $x\le  k<x+1$}
&=\frac{1}{n!\cos\theta}\sum_{k=1}^{n-1}vol\left(\left\{(x_1,\dots,x_{n-1})\left|
\begin{array}{l}
a\le x_1,\dots,x_{n-1}<1,  \\
k-1<\sum_{i=1}^{n-1} x_i\le k-a
\end{array}
\right.
\right\}\right).
\end{align*}
Putting $ x_i=a+t_i(1-a)$, we see that the summand is
\begin{align*}
&\int_{a\le x_1,\dots,x_{n-1}<1,  \atop
k-1<\sum_{i=1}^{n-1} x_i\le k-a} dx_1\dots dx_{n-1}
\\
=&(1-a)^{n-1}\int_{0\le t_i<1,\atop (k-na)/(1-a)-1\le \sum_{i=1}^{n-1}t_i\le (k-na)/(1-a)}
dt_1\dots dt_{n-1}
\\
=&(1-a)^{n-1}\{  U_{n-1}(\frac{k-na}{1-a}) -  U_{n-1}(\frac{k-na}{1-a}-1)  \}.
\end{align*}
Hence, by noting  $vol(\hat{\mathfrak{D}}_n)=1/(n!\cos\theta)$ (\cite{K1}), 
$V^a:=vol({D^a}\cap{\hat{\mathfrak{D}}_n})/vol(\hat{\mathfrak{D}}_n)$ is equal to
\begin{align*}
&
(1-a)^{n-1}\sum_{k=1}^{n-1}\left\{ U_{n-1}(\frac{k-na}{1-a}) -  U_{n-1}(\frac{k-na}{1-a}-1 )\right\}
\\
=&
\frac{(1-a)^{n-1}}{(n-1)!}
\sum_{k=1}^{n-1}   \left\{     \sum_{i=0\atop i< \frac{k-na}{1-a}}^{n-1} (-1)^i{\binom{n-1}i}
\left(\frac{k-na}{1-a}-i \right)^{n-1} \right.
\\
&\hspace{3.5cm}
-\left.   \sum_{i=0\atop i< \frac{k-na}{1-a}-1}^{n-1} (-1)^i{\binom{n-1}i}
\left(\frac{k-na}{1-a}-1-i \right)^{n-1} \right\}
\\
=&
\frac{(1-a)^{n-1}}{(n-1)!}
\sum_{k=1}^{n-1}   \left\{     \sum_{i=0\atop i< \frac{k-na}{1-a}}^{n-1} (-1)^i{\binom{n-1}i}
\left(\frac{k-na}{1-a}-i\right)^{n-1} \right.
\\
&\hspace{3.5cm}
+\left.   \sum_{j=1\atop j< \frac{k-na}{1-a}}^n (-1)^j{\binom{n-1}{j-1}}
\left(\frac{k-na}{1-a}-j\right)^{n-1} \right\}
\\
=&
\frac{(1-a)^{n-1}}{(n-1)!}   \left\{
\sum_{k=1\atop k> na}^{n-1}   \left(\frac{k-na}{1-a}\right)^{\hspace{-1mm}n-1} \right. 
\hspace{-3mm}
+\sum_{k=1}^{n-1}       \sum_{i=1\atop i< \frac{k-na}{1-a}}^{n-1} (-1)^i{\binom{n-1}i}
\left(\frac{k-na}{1-a}-i\right)^{\hspace{-1mm}n-1} 
\\
&\hspace{4cm}
+\left.\sum_{k=1}^{n-1}   \sum_{j=1\atop j< \frac{k-na}{1-a}}^{n-1} (-1)^j{\binom{n-1}{j-1}}
\left(\frac{k-na}{1-a}-j\right)^{n-1} \right\},
\intertext{
using ${\binom{n-1}i}+{\binom{n-1}{i-1}}=\binom ni$
}
=&
\frac{(1-a)^{n-1}}{(n-1)!}  \left \{
\sum_{k=1\atop k> na}^{n-1}   \left(\frac{k-na}{1-a}\right)^{\hspace{-1mm}n-1} \right.
+\sum_{k=1}^{n-1}  \left.      \sum_{i=1\atop i<\frac{k-na}{1-a}}^{n-1} (-1)^i{\binom{n}i} 
\left(\frac{k-na}{1-a}-i\right)^{\hspace{-1mm}n-1} 
 \right\}
 \\
=&
\frac{1}{(n-1)!}   \left\{
\sum_{k=1\atop k> na}^{n-1}   ({k-na})^{n-1} 
+      \sum_{i,k=1\atop i< \frac{k-na}{1-a}}^{n-1} (-1)^i{\binom{n}i}
({k-na}-i(1-a))^{n-1} 
 \right\} ,
 \\
 \intertext{replacing $i,k$ by $n-i,n-k$ in the second term}
 =&
 \frac{1}{(n-1)!}   \left\{
\sum_{k=1\atop k> na}^{n-1}   ({k-na})^{n-1} 
+      \sum_{k,i=1\atop (i-k)/i>a}^{n-1}  (-1)^{n+i}{\binom{n}i}
(i(1-a)-k)^{n-1} 
 \right\} 
 \\
  =&
 \frac{1}{(n-1)!}   \left\{
\sum_{k=1\atop k> na}^{n-1}   ({k-na})^{n-1} 
+      \sum_{1\le i\le n-1\atop i-n+1\le K \le i-1,K>ia}  (-1)^{n+i}{\binom{n}i}
(K-ia)^{n-1} 
 \right\} 
  \\
  =&
 \frac{1}{(n-1)!}   \left\{
\sum_{k=1\atop k> na}^{n-1}   ({k-na})^{n-1} 
+      \sum_{2\le i\le n-1\atop \max{(ia,i-n+1)}\le K \le i-1}  (-1)^{n+i}{\binom{n}{i}}
(K-ia)^{n-1} 
 \right\} 
  \\
  =&
 \frac{1}{(n-1)!}      \sum_{2\le i\le n\atop \max{(ia,i-n+1)}\le k \le i-1}  (-1)^{n+i}{\binom{n}{i}}
(k-ia)^{n-1} .
\end{align*}
Hence we have completed the proof.
\qed

\noindent
We note that $\left. \frac{dV^a}{da}\right|_{a=0}=-n$.
\begin{thm}
For $0\le a <1$, put
$$
E_a:=\{(x_1,\dots,x_n)\in[0,1)^n \mid x_n \le a\}.
$$
Then the right-hand side of \eqref{eq5} for $E_a$ is equal to
$$
\frac{1}{(n-1)!}\sum_{0\le i\le n\atop \max(ia,1)\le k \le n-1}(-1)^i\binom{n}{i}(k-ia)^{n-1}=V^{(1-a)}.
$$
\end{thm}
A proof is similar to the previous theorem, using the fact that
$$
\sum_{i=0}^n (-1)^{i}P(i)\binom{n}{i}=c_n(-1)^n n!
$$ for every polynomial $P(x)=c_nx^n +c_{n-1}x^{n-1}+\dots+c_0$ of degree at most $n$.

What is  the volume $vol(D\cap\hat{\mathfrak{D}}_n)$ for
$$
D:=\{(x_1,\dots,x_n)\in[0,1)^n\mid x_m\le a\}\,?
$$



\begin{thebibliography}{99}
%
\bibitem{Fe}
\sc{W. Feller}:\textit{An introduction to probability theory and its applications},
vol. 2, J. Wiley, New York, 1966.
%
%
\bibitem{K1}
\sc{Y. Kitaoka}:\textit{Statistical distribution of roots of a polynomial modulo primes II},
To appear U.D.T.
%
\bibitem{K2}
\sc{Y. Kitaoka}:\textit{Notes on the distribution of  roots modulo  a prime of a polynomial},
To appear U.D.T.
%
\end{thebibliography}
\end{document}